\magnification 1250
\topskip = 20 pt
\vsize=19.5cm
\hsize=13cm

\openup 1\jot

\headline{\ifnum\pageno=1 \hbox{}
    \else \ifnum\pageno<10 \hbox{\hskip 8truecm 801-0\folio\hfil}
          \else \hbox{\hskip 8truecm 801-\folio\hfil} \fi \fi}
\nopagenumbers

\noindent\hbox to\hsize{S\'eminaire BOURBAKI \hfil Novembre 1996} \par
\noindent 49\`eme ann\'ee, 1996-97, n$^{\rm o}$ 821 \par
\vskip 1.7truecm

\centerline{\bf PRODUCT FORMULAS FOR MODULAR FORMS ON $O(2,n)$}\par
\centerline{\bf [after R. Borcherds]} \par

\medskip
\centerline {by {\bf Maxim KONTSEVICH}}
\vskip 1.7truecm

\hyphenation{Lorent-zian}
\def\R{{\bf R}}
\def\Z{{\bf Z}}
\def\Q{{\bf Q}}
\def\C{{\bf C}}
\def\G{{\Gamma}}
\def\g{{\bf g}}
\def\H{{\cal H}}
\def\T{{\Theta}}
\def\L{{\Lambda}}
\def\l{{\lambda}}
\def\tl{{\tilde{\l}}}
\def\ltwo{{ c\bigl(- (\tl,\tl)/ 2\bigr) }}
\def\TL{{\tilde{\L}}}
\def\tp{{\tilde{p}}}
\def\ra{{\longrightarrow}}
\def\bs{{\backslash}}
\def\eps{{\epsilon}}

\noindent {\bf 1. INTRODUCTION }
\vskip 0.5truecm
\par \noindent {\bf 1.1. Product formulas }
\vskip 0.5truecm

A few years ago, R.~Borcherds found a remarkable multiplicative correspondence
 between classical modular forms with poles at cusps and meromorphic 
 modular forms on complex varieties 
$SO(n)\times SO(2)\bs SO(n,2)/\G\,,$ where $\G$ is an arithmetic 
subgroup in the real Lie group $SO(n,2)$. 
He was motivated by generalized Kac-Moody algebras, the Monster group and 
vertex operator algebras. The first proof of his 
formulas in completely classical terms (see [3])
 was rather indirect and complicated. 

In 1995 physicists J.~Harvey and G.~Moore wrote a paper on string duality where
 they found a new approach to Borcherds' identities (see [12]). They used divergent integrals, 
 which  look formally like integrals 
 in the classical theta correspondence in the theory of automorphic forms.
  R.~Borcherds recently wrote a preprint (see [5]) where he generalized his earlier results
 using the idea of Harvey and Moore. 
  My exposition is based mainly on this new preprint. 

Here is one of Borcherds' theorems:

\proclaim Theorem. Let $\L$ be an even unimodular lattice of signature $(s+1,1)$  
where $s=8,16,\dots$ and $v_0\in \L\otimes \R$ be a generic  vector of negative norm. 
Let $F=\sum_{-n_0}^{\infty} c(n) q^n\in \Z((q))$ be a meromorphic modular form of 
weight $(-s/2)$ for the group $SL(2,\Z)$ with poles only at the cusp.  
 Then there is a unique vector $\rho \in \L$ such that the function defined  
  for $v\in \L\otimes \C$ close to $it\,v_0,\,\,t\gg 1$, by the formula 
 $$
\Psi(v)=e^{2\pi i (\rho,v)}\prod_{\gamma\in L,\,\,\,(\gamma,v_0)>0} 
\left(1-e^{2\pi i (\gamma,v)}\right)^{c\left(-(\gamma,\gamma)/2\right)}$$
can be analytically continued to a meromorphic modular form of weight $c(0)/2$ 
for the group $O(s+2,2;\Z)^+$. In particular, the analytic continuation of $\Psi$ 
satisfies the equation
 $\Psi(2v/(v,v))=\pm\left((v,v)/2\right)^{c(0)/2}\Psi(v)$ .\par   

In my exposition I will formulate  results only in examples. One reason for 
this  is  that 
what is now known  is still far from the complete generality. Another reason is that 
I want to avoid heavy notations in order not to obscure the  logic of the  construction.
 
\vskip 0.5truecm
\par \noindent {\bf 1.2. An elementary example of a product formula }

\vskip 0.5truecm

Product formulas can be considered as statements about formal power  series of 
algebro-geometric origin. The general proof uses 
 analysis: integrals, infinite series and non-holomorphic functions. Here I 
will show a purely algebraic proof of a simple product formula. Both this formula 
and the proof are not new. They were discovered by D.~Zagier many years ago. Analogous 
formulas can be found in [11].
 
We fix  notations:
 $\H=\{\tau\in \C| \,\,Im(\tau)>0\}$ denotes the upper-half plane and
  $j:\H\ra\C$ is the classical elliptic invariant which identifies the quotient space 
  $\H/SL(2,\Z)$ with $ {\cal M}_1\simeq \C$, the coarse moduli space of complex elliptic 
  curves. We will compactify it to $\overline{\cal M}_1\simeq\C P^1$. 
  Function $q={\rm{exp}}\,(2\pi i \tau)$ can be considered
 as a holomorphic coordinate at a neighborhood of point $j=\infty$. 
  We expand the meromorphic function $j$ on $\C P^1$ in the coordinate $q$:
  $$j(q)=\sum_{n=-1}^{\infty} c(n) q^n=q^{-1}+744+\sum_{n=1}^{\infty} c(n) q^n\,\,\,\,.$$
\proclaim Theorem. For $0<|p|,|q|\ll 1$ one has
   the equality \hfil\break
\centerline{$ j(p)-j(q)=(p^{-1}-q^{-1})\prod_{k,l=1}^{\infty}(1-p^k q^l)^{c(kl)}$.} \par

From this equality follows an infinite sequence of algebraic identities between integer 
numbers $c(k),\,\,k\ge 1$. The first non-trivial identity is
  $$c(4)=c(3)+{c(1)^2-c(1)\over 2},\,\,\,\,\,\,
  20245856256=864299970+{196884^2-196884\over 2}$$

For each integer $n\ge 1$ there is an algebraic curve $C_n\subset \C\times \C$, the graph 
of the Hecke correspondence. In coordinates $(q_1,q_2)$ at the neighborhood
 of the point $(\infty,\infty)\in \C P^1 \times \C P^1$ this curve is again algebraic and 
 its branches are given by equations $q_1^{d_1}=q_2^{d_2}$ where $d_1d_2=n,\,\,d_1,d_2\ge 1$. 
 The Hecke operator $T_n$ is defined in the usual way using the correspondence $C_n$. It acts
 on meromorphic functions on $\C P^1$, on meromorphic 1-forms (=modular forms of weight 2), 
 etc. 

The main object will be a meromorphic differential 2-form on $\C P^1 \times \C P^1$ 
 $$\Omega=d_{along \,\,\,j_1} d_{along\,\,\, j_2}
  \left({\rm{log}}\,(j_1-j_2)\right)={1\over (j_1-j_2)^2} dj_1\wedge d j_2$$
where $(j_1,j_2)$ are coordinates on $\C P^1 \times \C P^1$. 

Let us write $j(q_1)-j(q_2)$ as an infinite product
 $(q_1^{-1}-q_2^{-1})\times \prod_{k,l=1}^{\infty} (1-q_1^kq_2^l)^{b(k,l)}$. We want to 
 prove that $b(k,l)=c(kl)$. Using symbols $b(k,l)$ we can write an explicit formula for the 
 $q$-expansion of $\Omega$:
$$\Omega={1\over (q_1-q_2)^2} dq_1\wedge d q_2 +\sum_{k,l=1}^{\infty} b(k,l) 
{kl\, q_1^{k-1} q_2^{l-1}\over (1-q_1^kq_2^l)^2} dq_1\wedge dq_2\,\,\,\,.$$

Now we  calculate the following double residue for $N,M\ge 1$:
 $$Res_{q_2=0}\, Res_{q_1=0}\left(j(q_1) j(q_2)\left[T^{(1)}_N\circ 
 T^{(2)}_M(\Omega)\right]\right)\,\,\,\,.$$

 This expression is equal to $0$ because $Res_{q_1=0}(\dots)$ is a meromorphic 
 1-form on $\C P^1$ with pole only at $j_2=\infty$. 
  We replace $\Omega$ by the sum as above. In the first term 
$$Res_{q_2=0} \,Res_{q_1=0}\left(j(q_1) j(q_2)\left[T_N^{(1)}\circ T_M^{(2)}\left( 
{1\over (q_1-q_2)^2} dq_1\wedge d q_2\right)\right] \right)$$  we can substitute $q_1^{-1}$ 
for $j(q_1)$ because of the regularity at $q_1=0$ of all other factors for generic $q_2$. 
Thus the first term can be expressed linearly in numbers $c(n)$. The second term is equal to
 $$Res_{q_2=0}\, Res_{q_1=0}\left(j(q_1) j(q_2) \sum_{k,l=1}^{\infty} (\dots)\right)=Res_{q_2=0}\,
  Res_{q_1=0}\left(q_1^{-1} q_2^{-1}\sum_{k,l=1}^{\infty} (\dots)\right)$$
 because of the regularity at zero of the double sum. This term can be expressed linearly in 
  numbers $b(k,l)$.  I  leave to the reader the rest of the calculation.

\vskip 1.2truecm
\noindent {\bf 2. STANDARD FACTS ABOUT AUTOMORPHIC FORMS}
\vskip 0.5truecm

\par \noindent {\bf 2.1. Definition of automorphic forms}
\vskip 0.5truecm

Let $G$ be a connected unimodular Lie group, $K$  a maximal compact subgroup, and $\G$
 a discrete subgroup of $G$ of finite covolume. Let us fix a homomorphism
 $\chi:{\cal Z} (U(\g))\ra \C$ from the center of the universal enveloping 
 algebra of the Lie algebra $\g:=Lie(G)$  to $\C$. Automorphic forms are 
 complex-valued $C^{\infty}$-functions
 on $G/\G$ which are $K$-finite and  annulated by a finite power of the ideal
  $Ker(\chi)$. A more general definition is obtained if one considers
  not just functions but sections of a local system associated with a 
  finite-dimensional representation $\rho:\G\ra GL(N,\C)$. 
  Any automorphic form is automatically real analytic because it satisfies
  an elliptic differential equation with real analytic coefficients. 

  Usually people consider functions satisfying certain growth conditions at 
   cusps, i.e. they consider $l^2$-integrable functions  or functions with 
    polynomial growth at infinity. In the classical case of $G=SL(2,\R)$ and
    $\G=SL(2,\Z)$ automorphic forms  include (anti)-holomorphic  modular forms
    of weights $k=4,6,...$ and Maass wave forms (eigenfunctions of the Laplace 
    operator on $\H/\G=K\bs G/\G$). The standard growth condition
     can  be formulated for holomorphic modular forms in terms of $q$-expansion as 
     the absence of terms $c(k) q^k$ with $k<0$.
      One of reasons to ignore automorphic forms with  exponential growth at cusps is 
      that the algebra
      of Hecke operators acts ``freely'' on such forms.

\vskip 0.5truecm
\par \noindent {\bf 2.2. Theta correspondence}
\vskip 0.5truecm

Theta correspondence transforms automorphic forms from one Lie group $G_1$ to another Lie 
group $G_2$, where $G_1\times G_2$ is a subgroup of the symplectic linear group (see [14]). 
The typical example is $G_1=Sp(V_1)$ and $G_2=SO(V_2)$ where 
$\left(V_1,(,)_1\right)$ is a symplectic real vector space and $\left(V_2,(,)_2\right)$ is a 
real vector space  with a non-degenerate symmetric bilinear form. The tensor product 
$V=V_1\otimes V_2$
     carries the natural symplectic structure $(,)_1\otimes (,)_2$. 

We denote by $W=W(V)$ the Hilbert space of the Weil  representation of the double covering 
$\widetilde{ Sp(V)}$ of the
     symplectic group $Sp(V)$. The space $W$ can be naturally identified with the space of 
     $l^2$-integrable functions on any Lagrangian subspace of $V$. Thus one can speak about 
     the nuclear space $W^{-\infty}$ consisting of distributions of  moderate growth. 
     Restricting the Weil representation of $\widetilde{ Sp(V)}$ to
 $\widetilde{ Sp(V_1)}\times \widetilde{ SO(V_2)}$ one gets a partially defined correspondence 
 between
 projective representations of $G_1$ and $G_2$. R.~Howe observed that this is
 a partial bijection in many cases. 
 
Let  
 $\Lambda_1\subset V_1$ and $\Lambda_2\subset V_2$ be integral lattices. Then 
 $\Lambda:=\Lambda_1\otimes \Lambda_2$ is an integral lattice in $V$. Denote by 
 $\Gamma_1,\,\Gamma_2,\,\Gamma$ arithmetic subgroups of $G_1,\,G_2,\,G=Sp(V)$ consisting 
 of automorphisms of these lattices. The space of invariants 
 $\left(W^{-\infty}\right)^{\Gamma}$ is finite-dimensional and consists of certain 
 theta functions. Theta correspondence is given by an integral operator from $G_1/\G_1$ 
 to $G_2/\G_2$ with the kernel equal to a theta function. In the next section we will
 consider an important example.

\vskip 0.5truecm
\par \noindent {\bf 2.3. Siegel theta function}
\vskip 0.5truecm

Let $(V_1,\Lambda_1)$ be $(\R^2,\Z^2)$ with the standard symplectic form and $\Lambda_2$ 
be an even unimodular lattice of signature $(n_+,n_-)$. We denote by 
$Gr$ the set of orthogonal decompositions of $V_2:=\Lambda_2\otimes \R$
 into  the  sum $V_+\oplus V_-$ of positive definite and negative definite
  subspaces. 
  $Gr$ can be considered as 
  an open subset of the Grassmanian of $n_-$-dimensional subspaces in
   $V_2$. 
   If $p\in Gr$ is such a decomposition  we denote by
 $p_+,\,\,p_-$ projectors onto $V_+,\,\,V_-$ respectively.
  The Siegel theta function (see [16]) is the restriction of 
  the standard theta function for $Sp(V_1\otimes V_2)$ to  the symmetric subspace
   $\H \times Gr\subset
   Sp(V_1\otimes V_2)/U(n_+ + n_-)$. 
  The explicit formula for it is
 $$\T(\tau,p)=\sum_{\lambda\in \Lambda_2} {\rm{exp}}\,\left( 2\pi i \left( 
 {(p_+(\l),p_+(\l))\over 2}\tau+{ (p_-(\l),p_-(\l))\over 2}\overline{\tau} \right)\right)=$$
$$=\sum_{\lambda\in \Lambda_2} q ^{(\l ,\l )\over 2} |q|^{-(p_-(\l),p_-(\l))}\,\,\,\,.$$
 
 This function is invariant under the action of $\G_2=Aut(\Lambda_2)$ on $Gr$. It transforms 
 under the action of $\G_1=SL(2,\Z)$ on $\tau\in \H$ as 
 $$\T\left({a\tau +b\over c\tau +d},p\right)=\pm (c\tau+d)^{n_+/2} (c\overline{\tau}+d)^{n_-/2} 
 \T(\tau,p)\,\,\,\,.$$

If $F$ is a holomorphic modular form for $\G_1=SL(2,\Z)$ of weight ${n_--n_+\over 2}$, or a 
Maass form for $n_-=n_+$, then  the theta transform of $F$ is defined as
$$\Phi(p)=\int\limits_{\H/PSL(2,\Z)} \T(\tau,p) F(\tau) y^{n_-\over 2} {dx dy\over y^2}\,\,\,,$$
where $\tau=x+iy, \,\,x,y,\in \R$. This integral converges for parabolic $F$.

The image of theta transform satisfies differential equations. Namely, there is a homomorphism 
$\alpha: {\cal Z}(U(\g_2))\ra {\cal Z}(U(\g_1))$ such that for any vector $v$ in the Weil 
representation $W(V_1\otimes V_2)$ and  any $z\in {\cal Z}(U(\g_2))$ one has 
$z(v)=(\alpha(z))(v)$. If we apply it to the theta function  we get the formula for the 
annulator of $\Phi$ in ${\cal Z}(U(\g_2))$.

The Siegel theta function also appears in string theory where it is the partition function of 
the torus $\C/(\Z+ \tau \Z)$ in  the Narain model associated with
 the indefinite lattice $\Lambda_2$ and the orthogonal decomposition $p$.

\vskip 1.2truecm
\noindent {\bf 3. BORCHERDS-HARVEY-MOORE CONSTRUCTION}

\vskip 0.5truecm
\par \noindent {\bf 3.1. Classical modular forms with poles at cusps}
\vskip 0.5truecm

The main idea of Borcherds-Harvey-Moore construction is a formal application of theta 
correspondence to modular forms for arithmetic subgroups in $SL(2,\R)$ with at most 
exponential growth at the cusps.
   In holomorphic case and for $\G_1=SL(2,\Z)$ any such form is meromorphic on
    $\H/\Z\sqcup\{\infty\}$ 
   with only pole at the cusp. It  can be presented as
 $\Delta(\tau)^{-k} F_0(\tau)$ where $k\ge 0$, 
 $\Delta=q\prod_{n=1}^{\infty}(1-q^n)^{24}$, and $F_0$ is a parabolic 
 holomorphic modular form.  Unlike in the classical theory, 
 holomorphic modular forms with poles at cusps can have negative weights.
 In the case of Maass forms
 for any $\lambda\in \C$ there is an infinite-dimensional vector space of solutions of the 
 equation $\Delta F=\lambda F$ on $\H/\G_1$ with exponential growth at the cusp.

  There are also other automorphic forms like
$E_2'(\tau):=1-24 \sum_{n=1}^{\infty} {n\, q^n\over 1-q^n}-{3\over \pi y}$, real analytic 
Eisenstein series, Siegel theta functions, etc.

  All these forms have the following common property: 
   there exists $M\ge 0$
  such that for any $N\ge 0$ the form can be expanded in a neighborhood of the cusp as 
$$o(y^{-N})+\sum_{m:|m|<M} e^{2\pi (imx +|m| y)} \left( \sum_{j\in \,\,finite 
\,\,\,\,set} c_{m,j} ({\rm{log}}\,(y))^{k_{m,j}} y^{{\sigma}_{m,j}}+\eps_m(y)\right)$$
where $k_{m,j}\in \Z_{\ge 0}$, $m\in \Z$, $c_{m,j},\,\sigma_{m,j}\in \C$ and $\eps_m(y)=o(y^{-N})$ 
depends on $y$ only. 

\vskip 0.5truecm
\par \noindent {\bf 3.2. Regularization of divergent integrals}
\vskip 0.5truecm

Let us assume that $F$ is a holomorphic modular form of weight $(n_- -n_+)/2$ with poles at cusps. 
We want to make  sense of the divergent integral
$$\Phi(p)=\int\limits_{\H/PSL(2,\Z)} \T(\tau,p) F(\tau)\, y^{n_-\over2} {dx dy\over y^2}\,\,\,\,.$$
 
 After the expansion of $\T(\tau,p)$ and $F(\tau)$
  at $y=\rm{Im}\, \tau\ra+\infty$ 
  there are only finitely many divergent terms of the form
$$\int\limits_{x\in [0,1],\,\,\,y\ge const} {\rm{exp}}\,(2\pi i m x+2\pi |m|y 
-Ly)\,y^{n_-/2-2} 
\,\,\,dx dy\,\,\,\,,$$
where $L$ is a non-negative real-valued function on $Gr$.

If $m\ne 0$ then we define the regularized value of this integral to be  $0$.
 It is the  natural choice if we perform the integration along the variable $x$ first. If $m=0$ 
 and $L>0$ then the integral is absolutely convergent. The trouble arises only when $m=L=0$. 
 In this case we can subtract the divergent term $y^k$ from the
 integral in domain $y>y_0$. The result will be a function of $y_0$.

In the  more general case for non-holomorphic forms there are finitely many 
divergent components of the product 
$F(\tau)\T(\tau,p)$ with  frequency $m=0$ along coordinate $x$. All these 
components are of type $\bigl({\rm{log}}\,(y)\bigr)^k y^{\sigma}$. One 
can multiply $F$ by $y^{-s}$, or (better) by a real analytic Eisenstein series. Then we can 
assign a value to the integral for $ Re(s)\gg 0$ and continue it to a meromorphic function 
for all $s\in\C$. The regularized integral can be defined
 as the constant term of the Laurent expansion at $s=0$.

\vskip 0.5truecm
\par \noindent {\bf 3.3. Automorphic forms with singularities at locally homogeneous 
submanifolds}\vskip 0.5truecm

 Let us see what kind of divergences our integral has for holomorphic form $F=\sum_n c(n) q^n$.  
First of all, if the $c(0)\ne 0$  then the constant term in the series for $\T$ corresponding 
to the vector $\lambda=0$ produces troubles.  This problem is independent of the point $p$ in 
the Grassmanian, and we can resolve it in one way or another. The result is that we still can 
define $F$  modulo an additive constant.

Other divergent terms appear when there is a \it non-zero \rm lattice vector 
$\lambda\in\Lambda_2 $ such that $\lambda$ belongs to the the positive subspace $V_+$ and has 
a special length. Namely, a term $q^{-(\lambda,\lambda)/2}$ should be present in the $q$-expansion 
of $F$.
  
Thus we see that the singular set of $\Phi$ in $X=Gr/\G_2$ consists of a finite union of 
certain totally geodesic submanifolds of type $X'=K'\backslash G'/\G'$. The same fact holds for
 non-holomorphic modular forms $\Phi$ admitting an asymptotic expansion at infinity as
  in 3.1. Also, one can
 write explicitly the types of singularities of $\Phi$, i.e. functions $\Phi'$ defined at a 
 neighborhood of $X'$ such that $\Phi-\Phi'$ can be continued to a real-analytic function. 
 These functions $\Phi'$ are finite linear combinations 
 of functions $x\mapsto({\rm{log}}\,(dist(x,X'))^k (dist(x,X'))^{\sigma}$ 
 where $dist(x,X')$ 
 is the distance between $x$ and $X'$ in a natural metric.

 The function $\Phi$ on the domain of definition satisfies differential equations.  If $c(0)=0$ 
 then these will be homogeneous linear differential  equations $z(\Phi)=0$ for some 
 $z\in{\cal Z}U(\g_2)$ (see the end of 2.3). The  divergent term $c(0)q^0$ produces certain 
 universal r.h.s. for these equations.

In the case $n_+=1$ submanifolds $X'$ are locally hyperplanes in the hyperbolic space.
 In the case $n_+=2$ they are complex hypersurfaces. The same is true for $G_2=SU(N,1)$. 
 If $G_2=SP(2g,\R)$ and  we consider $G_1=PSL(2,\R)$ as the orthogonal group $SO(2,1)$) then 
 subvarieties $X'$ are complex subvarieties of codimension $g$.

If $(G_2,\G_2)=(SL(2,\R),SL(2,\Z))$ then the forms $\Phi$ on $\overline{\cal M}_1$ have 
singularities at 
 Heegner points, i.e. moduli of elliptic curves with complex multiplication, or equivalently, 
 points with the coordinate $\tau$ in an imaginary  quadratic field,
 $\tau=x+iy$, where $x\in \Q$ and $y^2 \in \Q$.

\vskip 0.5truecm
\par \noindent {\bf 3.4. Fourier expansions at cusps}
\vskip 0.5truecm

In theta correspondence one can write an expansion of $\Phi$ at cusps via Fourier coefficients 
of $F$. The usual trick (Rankin-Selberg method) consists in replacing of the integral over the 
fundamental domain of $SL(2,\Z)$ by an integral over the fundamental domain of $\Z$ in $\H$. In 
the case of divergent integrals one should be cautious when interchanging infinite sums and 
integrals.

Let $\l_0\in\Lambda_2$ be a primitive 
null-vector, $(\l_0,\l_0)=0$. Such vector always exists
 for undefinite lattices of sufficiently large rank, including all
  even unimodular lattices.
   We define smaller lattice $\tilde{\L}$ as $\l_0^{\perp}/\Z \l_0$. 
It is easy to see that any orthogonal  decomposition $p$ of $V_2=\Lambda_2\otimes \R$ defines an 
orthogonal decomposition $\tilde p$ of $\tilde{V}= \tilde{\L}\otimes \R$ using natural
isomorphism between
 $\tilde{V}$ and $\left(\R \,p_+(\l_0)+\R \,p_-(\l_0)\right)^{\perp}$.

Using Poisson summation formula in direction $\Z \l_0$ one can rewrite $\T(\tau,p)$ as  
$$\T(\tau,p)={1\over \sqrt{2 y (p_+(\l_0),p_+(\l_0))} } \sum_{\l'\in \L_2/\Z \l_0} 
\sum_{l\in \Z} {\rm exp}\,
(\dots)$$
where ${\rm exp}\,
(\dots)$ is the exponent of an explicit algebraic expression.
 
Let us choose an additional lattice vector $\l_1$ such that $(\l_0,\l_1)=1$. 
Then  we can  embed $\tilde{\L}$ in $\L_2$ as $(\Z\l_0+\Z\l_1)^{\perp}$. Moreover, 
we can parametrize $\L_2/\Z \l_0$ by $\tilde{\L}\times \{k \,\l_1|\,k\in \Z\}$.
  The total sum becomes a sum over $(k,l)\in \Z^2$ of certain twisted theta series 
  associated with the lattice $\tilde{\L}$ and the orthogonal decomposition $\tilde p$. 
  The total formula is quite cumbersome.  

\proclaim Main Identity. Let $F$ be a bounded measurable
 function on $\H$ satisfying the equation 
$$F\left({a\tau+b\over c\tau+d}\right)=
(c\tau+d)^{n_- -n_+\over 2} F(\tau)\,\,\,.$$
Then the following identity holds:
$$\int\limits_{\H/SL(2,\Z)} \T(\tau,p) F(\tau) y^{n_-\over 2}
 {dx dy\over y^2}\,\, =\,\,
 {1\over \sqrt{2\,\epsilon} } \int\limits_{\H/SL(2,\Z)}  
 \T(\tau,\tilde{p}) F(\tau) y^{n_- -1\over 2} {dx dy\over y^2}\,\,+$$
$$+{2\over \sqrt{2\,\epsilon}} \sum_{n>0} \,\, \int\limits_{\H/\Z} F(\tau) 
\,\,{\rm{exp}}\,\left(-{\pi n^2\over 2 y\,\epsilon}\right) \sum_{\tl\in 
\TL} e^{2\pi i n(\tilde{\l},\mu)} q^{(\tilde{\l},\tilde{\l})/2} 
|q|^{-(\tilde{p}_-(\tilde{\l}),\tilde{p}_-(\tilde{\l}))} y^{n_- -1\over 2} 
{dx dy\over y^2}\,\,.$$\par

Here $\epsilon=(p_+(\l_0),p_+(\l_0))>0$ and 
$$\mu=\left(\l_1+(p_+(\l_0)-p_-(\l_0))/2\epsilon\right) \,mod\,\R\l_0\in 
{\l_0}^{\perp}/\R\l_0=\tilde{\L}\otimes \R\,\,\,.$$
 Correspondence $p\mapsto(\tp,\eps,\mu)$ can be considered as a local
  parametrization of $Gr$.

This identity  follows 
from the formula for $\T(\tau,p)$ as a sum over pairs of integers $(k,l)$. 
The first term comes from the term corresponding to $k=l=0$. 
Any other pair of integers $(k,l)$ can be presented as $(nc,nd)$
 where $n>0$ and $c$ and $d$ are coprime. We identify in the usual way
 $SL(2,\Z)/\left\{\left(\matrix 
 {1 & n\cr 0 & 1}\right)\Big |\,n\in \Z\right\}$ 
 with the set of primitive vectors $(c,d)$ in $\Z^2$ and rewrite the sum 
 over non-zero pairs $(k,l)$ as the sum over $n>0$ and over the set of 
 copies in $\H/\Z=\{(x+iy|\,0<x\le 1,\,\,\,0<y<+\infty\}$ of the classical fundamental domain 
 of $SL(2,\Z)$ in $\H$.

Now let us see what happens for functions $F$ which admit an asymptotic expansion at the cusp 
as in 3.1. The divergence of the sum above as $y\ra\infty$ is of the same form as one of the 
original integral for $\Phi$. In the integral over $\H/SL(2,\Z)$ one might expect a priori 
divergences near cusps on $\Q\subset\R$. In fact, the function $F$ has exponential growth at 
these points. Nevertheless the integral near $\R$ is absolutely convergent because of the 
factor ${\rm{exp}}\,(-\pi/2 \epsilon y )$ which makes the total integrand small enough, as 
$\epsilon\ll 1$. Thus the exchange of the order of the sum and of the integral is justified 
for small
 $\epsilon=(p_+(\l_0),p_+(\l_0))$.

One can calculate explicitly integrals over $\H/\Z$ corresponding to individual
 terms in $q$-expansion of $F$. It reduces to classical integrals for Bessel functions 
$$\int_{y>0} {\rm{exp}}\,(-\beta/y-\alpha y) y^{\nu-1} dy= 2(\beta/\alpha)^{\nu/2} K_{\nu}
(2\sqrt{\alpha\beta})\,\,\,\,.$$ 

\vskip 0.5truecm
\par \noindent {\bf 3.5. Hyperbolic case}
\vskip 0.5truecm

Let us consider the case of hyperbolic even unimodular lattices,
 $n_-=1$ and $n_+>1$. As in the previous section, we pick a primitive 
 null-vector $\l_0\in \L_2$. The additional vector $\l_1$ such that
  $(\l_0,\l_1)=1$ is chosen  now among null-vectors. Lattice 
   $\TL\simeq (\Z\l_0+\Z\l_1)^{\perp}$ 
   is considered as a sublattice of $\L$. Thus we have a
   decompositon
   $\L_2=\Z\l_0+\Z\l_1+\TL$.
   
   We identify the space $Gr$ with the (half of the) hyperboloid
    $$H=\{v\in \L_2\otimes \R| \,\,(v,v)=-1,\,\,\,(v,\l_0)>0\}\,\,\,.$$
    Projector $p_-$ corresponding to $v$ is the orthogonal projector
     to $1$-dimensional space $\R \,v$.
     In terms of parameters $(\eps,\mu)\in \R_+\times
     (\TL\otimes \R)$ from the previous section we have
     $$v=v(\eps,\mu)={- {1\over \sqrt{\eps}} -\sqrt{\eps}(\mu,\mu)
     \over 2} \l_0+\sqrt{\eps}\l_1+\sqrt{\eps}\mu\,\,\,.$$
     Parameter $\tp$ does not vary because $\TL$ is positive definite and 
     $\tilde{Gr}$ is a one point set.
     
     Let $F=\sum_n c(n) q^n\in \C((q))$ be a holomorphic modular form of weight
      $1-n_+\over 2$. 
      
      \proclaim Theorem. Theta transform $\Phi$ of $F$ is locally 
      the restriction 
      of a continuous piecewise linear function on $\TL\otimes \R$ to the
       hyperboloid $H$.\par
       
      The rest of this section is devoted to the proof of this Theorem.
      Function $\Phi$ has the following expansion as $\eps\ra 0$:
     $$\Phi=\Phi(\eps,\mu)=
     {1\over \sqrt{2\,\epsilon} } \int\limits_{\H/SL(2,\Z)}  
 \T(\tau,\tilde{p}) F(\tau) {dx dy\over y^2}\,\,+$$
   $$ +{2\over \sqrt{2\,\epsilon}} \sum_{n>0} \,\, \int\limits_{\H/\Z} F(\tau) 
\,\,{\rm{exp}}\,\left(-{\pi n^2\over 2 y\,\epsilon}\right) \sum_{\tilde{\l}\in 
\tilde{\L}} e^{2\pi i n(\tilde{\l},\mu)} q^{(\tilde{\l},\tilde{\l})/2} 
{dx dy\over y^2}\,\,.$$  
  Notice that each integral in this formula can be
   unambigously regularized using rules from 3.2. We denote terms in this
    formula by $\Phi_{(1)}$ and $\Phi_{(2)}$.
   
   Term $\Phi_{(1)}$ is proportional to
   $(v,\l_0)^{-1}$ because 
  $(v,\l_0)=\sqrt{\eps}$.

      After the expansion of $F$ into the series we see that in
       the  term $\Phi_{(2)}$
      we have to calculate integrals
     $$\int\limits_{\H/\Z} q^{m+{(\tl,\tl)\over 2}} {\rm{exp}}\,\left(-{\pi 
     n^2\over 2 y\,\epsilon}\right) {dx dy\over y^2}\,\,.$$
     If $m+{(\tl,\tl)\over 2}=0$ then  this integral is equal to $2\eps/ \pi n^2$,
      otherwise it vanishes.
      Thus we see that the second term is 
      $$\Phi_{(2)}={4 \pi \sqrt{\eps}\over \sqrt{2}} \sum_{\tl\in \TL}
      \ltwo \sum_{n>0} {e^{2\pi i n(\tilde{\l},\mu)}\over n^2}\,\,\,.$$
      Every vector $\tl$ appears in this formula together with the opposite
       vector $-\tl$. It implies that we can replace in the formula from above
        exponent by the cosine. Now we use the formula
        $$\sum_{n>0} {{\rm cos}\,(2\pi n x)\over n^2}=\pi^2\left(x^2 +\alpha(x) x+{1\over 6}
        \right),\,\,\,x\in \R$$
        where $\alpha(x)$ is a locally constant function of $x$: $\,\alpha(x)=-2n-1$ for $n\le x<n+1$.

         Finally, we get a formula for $\Phi$ as a finite sum:
         $$\Phi(\eps,\mu)=\Phi_{(1)}+
         {\pi\over \sqrt 2}
         \sum_{\tl\in \TL} \ltwo\cdot  4{\sqrt{\eps}}
         \left\{ 
         (\tl,\mu)^2+\alpha((\tl,\mu)) (\tl, \mu)+{1\over 6}
         \right\}
         \,\,\,.$$
       
       Terms (locally) proportional to $\sqrt{\eps}(\tl,\mu)$ are restrictions
        of linear functions $v=v(\eps,\mu)\mapsto const(v,\tl)$. Terms proportional
         to $\sqrt{\eps}$ are restrictions of linear functions $v\mapsto const (v,\l_0)$.
         We claim that the rest is also the restriction
          of linear function (proportional to $v\mapsto (v,\l_1)$ in fact).

         Using the fact that $\Phi_{(1)}(v)=const \,(v,v)/(v,\l_0)$
         we see that 
          $$\Phi=\Phi(v)=({\rm piecewise\,\,\,linear\,\,\,function
          \,\,\,of\,\,\,}v)+{1\over (v,\l_0)}({\rm quadratic\,\,\,
          polynomial\,\,\,of\,\,\,}v)\,\,\,.$$
          
        Applying the following automorphism of $\L_2$:
        $$\l_0\mapsto\l_1,\,\,\l_1\mapsto \l_0,\,\,\tl\mapsto\tl\,\,\,{\rm for}
        \,\,\,\tl \in\TL$$
        we see that 
        $$\Phi(v)=({\rm piecewise\,\,\,linear\,\,\,function
          \,\,\,of\,\,\,}v)+{1\over (v,\l_1)}({\rm quadratic\,\,\,
          polynomial\,\,\,of\,\,\,}v)\,\,\,.$$
          
          Comparing two expressions for $\Phi(v)$ as above we conclude
           that $\Phi$ is a piecewise linear function of $v$.

      This finishes the proof of the Theorem of this section.
       R.~Borcherds calculated theta transform in a more general situation
        and obtained that $\Phi$ is the restriction of a  piecewise polynomial function on the
         hyperboloid $H$.

\vskip 0.5truecm
\par \noindent {\bf 3.6. Product formulas for meromorphic forms}
\vskip 0.5truecm

 Now we consider the case when $n_-=2$ and $F$ is a holomorphic modular form with  pole at the cusp. 
 
 As in the previous section we fix decomposition $\L_2=\Z\l_0+\Z\l_1+\TL$ where $\TL$ is now
  a hyperbolic lattice.
  The space $Gr$ is convenient to parametrize by vectors $v\in \TL\otimes \C$
   such that
  $( Im(v), Im(v))<0$. Projector $p$ in $\L_2\otimes \R$ and corresponding parameters $(\tp, \eps,\mu)$
   from the Section 3.4 are given by the following formulas:
   $$p_-(\L_2\otimes \R)=\R\cdot Re(u)+ \R\cdot Im(u)
   \subset \L_2\otimes \R \,\,\,\,{\rm  where}\,\,\,\, u=-{(v,v)\over 2}\l_0+\l_1+v,$$
    $$\tp_-(\TL\otimes \R)=\R\cdot Im(v),$$
   $$\eps={-1\over (Im(v),Im(v))}.$$
   $$\mu=-Re(v)\,\,\,.$$

The integral over $\H/ SL(2,\Z)$ in the first term of the formula for 
 the Fourier expansion at cusps (Section 3.4) was evaluated in the pervious
  section. The result is that the first term has the form
   $$\Phi_{(1)}=\Phi_{(1)}(v)=(W(v),Im(v))$$
   where $W(v)$ is a locally constant function on $Gr$ with values in $\TL\otimes\C$.
  
Now we consider $\Phi_{(2)}$, the sum of integrals over $\H/\Z$.

The contribution of terms with $\tilde{\l}=0$ and $n>0$ is a divergent sum $const + 2\,c(0) \sum_{n>0}1/n$. 
Nevertheless, if we use some regularization procedure, we obtain 
 $2\,c(0) \,\sum_{n>0} \eps^s/n^{2s+1}$ as $s\ra 0$. An easy calculation shows that
 the regularized value is $const+c(0) {\rm{log}}\,(\eps)$. 

The contribution of the term corresponding to $\tilde{\l}\ne 0$ and 
$n>0$ 
 is equal to
$$2\,\ltwo{1\over n} \,e^{2\pi i n(\tl,-Re(v))}\,
{\rm{exp}}\,\left(-2\pi n\,\sqrt{-(\tp_-(\tl ),\tp_-(\tl)\over \eps}\right)\,\,\,\,.$$
Here we use the classical formula $K_{-1/2}(z)=\sqrt{\pi/2z} \cdot {\rm exp}\,(-z)$. 

Elementary
 calculations show that $\sqrt{-\bigl(\tp_-(\tl ),\tp_-(\tl)\bigr)\over \eps}=|(\tl , Im(v))|$.
 The sum over $n$ and two opposite vectors $\pm\tl$ gives 
  $$-4\,\ltwo\,{\rm{log}}\,\,|1-{\rm{exp}}\,(2\pi i (\tl_+,v )|$$ 
  where $\tl_+$ is the one of two vectors $(\tl,-\tl)$ which has positive scalar product with $Im(v)$.
  
  The resulting formula for $\Phi$ is (up to an additive constant)
  $$(W(v),Im(v))+c(0) {\rm{log}}\,(\eps)-4
  \sum_{\tl:\,(\tl,Im(v))>0} \ltwo\,{\rm{log}}\,\,
  |1-{\rm{exp}}\,(2\pi i (\tl_+,v )|\,\,\,.$$
  Notice that all terms here except $c(0)\, {\rm{log}}\,(\eps)$ are locally sums of holomorphic
   and anti-holomorphic functions of $v$.

Let us  now assume that the coefficients $c(n)$ of $F$ are integers.
 Denote by $\cal L$ an equivariant complex line bundle over $Gr$ whose total space 
 is the complex cone 
 $$\{w\in \L_2\otimes \C|\,(w,w)=0,\,\,\bigl(Im(w),Im(w)\bigr)< 0\}\,\,\,.$$ Projection 
 ${\cal L}\ra Gr$ is $ w\mapsto u={w / (w,\l_0)}$. Bundle $\cal L$ carries
   invariant hermitean scalar product $\|w\|:=\sqrt{ -(w,{\overline w})/2}$.
  We claim that there exists a meromorphic section $\Psi$ of $({\cal L})^{\otimes (c(0)/2)}$
   such that 
   $${\rm{log}}\,\|\Psi\|=-\Phi/4\,\,
   \,.$$
 Locally, it follows from the expression for $\Phi$ from above and from the identity
 $$\eps={  |(w,\l_0)|^2\over \|w\|^2}\,\,\,.$$
Globally, we use the information about singularities from 3.3. 
In general, $\Psi$ is not 
 $\G_2$-equivariant and it gives a section  on $Gr/\G_2$ of ${\cal L}^{\otimes (c(0)/2)}$ twisted  
 with a unitary character of $\G_2/[\G_2,\G_2]$.  In this way one obtains a proof of
  the Theorem from 1.1.

 R.~Borcherds proved  more general product formulas for congruence 
 subgroup in $SL(2,\Z)$, non-unimodular lattices, and proposed to consider the case $G_2=SU(N,1)$. 
 He also developed the formalism for generalized theta functions associated with harmonic polynomials.

\vskip 1.2truecm 
\noindent {\bf 4. EXPLICIT EXAMPLES}
\vskip 0.5truecm

I will present only $3$ examples.

 The simplest example of the product formula for the group $SO(1,2)$ is completely trivial:
  $\eta(q)=q^{1/24}\prod_{n=1}^{\infty} (1-q^n)^1$
 where all the exponents $1$ are coefficients of $q^{n^2}$ of a form of weight $1/2$,
 namely of the theta function $1/2+\sum_{n=1}^{\infty}q^{n^2}$.

An example for $SO(2,2)$ is the product formula for $j(p)-j(q)$ in 1.2.

  The next  example  for the group $SO(3,2)$ is also very beautiful:
$$\sum_{m,n\in \Z} (-1)^{m+n} p^{m^2} q^{n^2} r^{mn}=\prod_{a+b+c>0}\left({1-p^a q^c r^b\over 
1+p^a q^c r^b}\right)^{f(ac-b^2)}$$
where $\sum f(n) q^n =1/(\sum_n (-1)^n q^{n^2})=1+2q+4q^2+8q^3+14q^4+\dots$.

\vskip 1.2truecm
\noindent {\bf 5. CONNECTIONS WITH OTHER PARTS OF MATHEMATICS}
\vskip 0.5truecm
\par \noindent {\bf 5.1. Generalized Kac-Moody algebras }
\vskip 0.5truecm

There are many nice examples of so-called generalized Kac-Moody algebras constructed by 
R.~Borcherds and later by V.~Gritsenko  and V.~Nikulin. These Lie super-algebras are graded by a lattice, 
and the generating function for the dimensions of homogeneous components is essentially an 
automorphic form for $O(n,2)$. The product formula  is the Weyl-Kac-Borcherds denominator 
identity. The Weyl group for these algebras is often an arithmetic group generated by 
reflections. Also the Monster group appears
 as an automorphism group.

\vskip 0.5truecm
\par \noindent {\bf 5.2. K3 surfaces, Mirror symmetry, string duality}
\vskip 0.5truecm

Let $\L$ be an even hyperbolic sublattice of (unique) even unimodular 
lattice $\L_{3,19}$ of
 signature $(3,19)$. We define ${\cal M}_{\L}$ as the moduli space of algebraic $K3$-surfaces 
 $X$ such that $\L\subset Pic(X)\subset H^2(X,\Z)\simeq \L_{3,19}$.
  By the classification theorem for $K3$-surfaces we 
  see that ${\cal M}_{\L}$
 is of the type where Borcherds' products are defined. 
 The  image of ${\cal M}_{\L}$ under  the period map is 
 $$\{w\in \L^{\perp}\otimes \C|\,\,(w,w)=0,\,\,(w,\overline{w})<0;
 \,\,\forall \l\in\L_{3,19}\cap 
 \L^{\perp}\cap w^{\perp} \,\,\,\,\,\,\,(\l,\l)\ne-2\, \}/\C^{\times}\,\,\,.$$

Thus Borcherds' results mean that  for some $\L$ 
 the standard line bundle over ${\cal M}_{\L}$ is a torsion 
element in $Pic({\cal M}_{\L})$. One of such lattices is  the
one-dimensional lattice
 $\L=\Z\l,\,\,(\l,\l)=2$.
  In general, product formulas produce linear relations between certain divisors in Shimura varieties.

Also one can consider the moduli space of Riemannian metrics on $4$-dimensional manifolds, 
which are hyperkaehler metrics on $K3$-surfaces. This space is locally modeled by 
$SO(19,3)/SO(19)\times SO(3)$. Borcherds' construction gives
 a certain real-analytic function on it. Presumably, it is related to the
 regularized determinant of the Laplace operator (see [15]).

Gritsenko-Nikulin, Jorgenson-Todorov  and Harvey-Moore (see [9,13,15]) made  conjectures about 
the relation between Kac-Moody algebras, $K3$-surfaces, Borcherds' product formulas
 and mirror symmetry. One expects 
that certain  numbers of curves of various genera on  a generic element of the family 
${\cal M}_{\L}$ coincide with exponents in a product formula
 associated with the dual family where $\L$ is equal to the lattice of transcedental cycles on 
 a generic element.  
 
In fact, J.~Harvey and G.~Moore tried to find  a conceptual construction of the generalized
 Kac-Moody algebra associated with the $K3$-surface $X$. Conjecturally,
 it is the direct sum of all cohomology groups of all moduli spaces of stable
 coherent sheaves on $X$. The Lie bracket is given by a correspondence 
 in the cube of the total (disconnected) moduli space. This correspondence
 is expected to consist of triples of sheaves from all possible short exact sequences. 

The idea of the Harvey-Moore integral arose 
 from a new duality in string theory relating elliptic curves on one manifold
 to numbers of all curves of all genera on another manifold. Thus the integral over the moduli 
 space of elliptic curves appeared. In some cases we get in the formula for $\Phi$ an infinite 
 sum of $3$-logarithm functions as in the usual mirror symmetry. R.~Dijkgraaf, G.~Moore,
  E.~Verlinde and H.~Verlinde proposed a Borcherds' type identity involving elliptic genera, see [7].

\vskip 0.5truecm
\par \noindent {\bf 5.3. Hyperbolic case and Donaldson invariants}
\vskip 0.5truecm

It is well-known that in Donaldson theory $4$-dimensional manifolds $X$ with the $b_+=1$ are 
very special. The Donaldson invariant is a piecewise polynomial function
 on the cone $\{x\,|\,\,x\cdot x<0\}$
  in the hyperbolic space $H^2(X,\R)$. R.~Borcherds observed that in certain cases 
 (like $\C P^2$ blown up at $9$ points) the Donaldson invariant coincides
 with one of the functions $\Phi$ given by the theta correspondence.

\vskip 1.5truecm
\centerline{\bf REFERENCES}
\vskip 0.5truecm
\item{[1]} R.~Borcherds - {\it Generalized Kac-Moody algebras}, J. of Algebra {\bf 115} (1988), 501-512.
\item{[2]} R.~Borcherds - {\it Monstrous moonshine and monstrous Lie superalgebras}, 
Invent. Math. {\bf 109} (1992), 405-444.
\item{[3]} R.~Borcherds - {\it Automorphic forms on $O_{s+2,2}$ and infinite products}, 
Invent. Math. {\bf 120} (1995), 161-213.
\item{[4]} R.~Borcherds - {\it The moduli space of Enriques surfaces and the fake monster 
Lie superalgebra}, Topology {\bf 35} no. 3 (1996), 699-710.
\item{[5]} R.~Borcherds - {\it Automorphic forms with singularities on Grassmanians}, 
e-print alg-geom/9609022 (1996).

\item{[6]} A.~Borel, H.~Jacquet - {\it Automorphic forms and automorphic representations}, 
Proc. Symp. Pure Math. {\bf 33} part 1 (1979), 189-202.
\item{[7]} R.~Dijkgraaf, G.~Moore, E.~Verlinde, H.~Verlinde - {\it Elliptic genera of 
symmetric products and second quantized strings}, e-print hep-th/9608096 (1996).
\item{[8]} V.~Gritsenko, V.~Nikulin - {\it Siegel automorphic form corrections of some 
Lorentzian Kac-Moody algebras  }, e-print alg-geom/9504006 (1995).

\item{[9]} V.~Gritsenko, V.~Nikulin - {\it K3 surfaces, Lorentzian Kac-Moody algebras and 
Mirror Symmetry}, Math. Res. Letters {\bf 3} (1996), 211-229.
\item{[10]} V.~Gritsenko, V.~Nikulin  - {\it The Igusa modular forms and ``the simplest'' 
Lorentzian Kac-Moody algebras}, e-print alg-geom/9603010 (1996).
\item{[11]} B.~Gross, D.~Zagier - {\it On singular moduli}, Jour. f\"ur die Reine und ang. 
Math. {\bf 355} (1985), 191-220.
\item{[12]} J.~Harvey, G.~Moore - {\it Algebras, BPS states and strings}, e-print hep-th/9510182; 
Nucl. Phys.  {\bf B463} (1996), 315.
\item{[13]} J.~Harvey, G.~Moore - {\it On the algebras of BPS states}, e-print hep-th/9609017
\item{[14]} R.~Howe - {\it $\theta$-series and invariant theory}, Proc. Symp. Pure Math. {\bf 33} 
part 1 (1979), 275-285.

\item{[15]} J.~Jorgenson, A.~Todorov - {\it An analytic discriminant for polarized K3 surfaces}, 
preprint (1994).

\item{[16]} S.~Rallis - {\it On a relation between $\widetilde{SL}_2$ cusp forms  and automorphic 
forms on orthogonal groups}, Proc. Symp. Pure Math. {\bf 33} part 1 (1979), 297-314.

\vskip 1truecm
\hskip 8truecm Maxim KONTSEVICH\par
\medskip
\hskip 8truecm I.H.E.S.\par
\hskip 8truecm 35, route de Chartres\par
\hskip 8truecm F--91440  Bures--sur--Yvette\par
\hskip 8truecm FRANCE

\smallskip
\hskip 8truecm E--mail~: maxim@ihes.fr\par

\end